\documentclass{primus}

\usepackage{amsmath,amssymb}
\usepackage{color}
\usepackage{graphicx}
\usepackage[hidelinks]{hyperref}
\usepackage{booktabs}

\graphicspath{{images/}}

\title{ArcheoLab: A hands-on course\\ to engage students in applied mathematics}

\author{
Denise Schmutz\\
	Faculty of Mathematics\\
	University of Vienna\\
	Oskar-Morgenstern-Platz 1\\
	A-1090 Vienna, Austria\\
        denise.schmutz@univie.ac.at
	\and
Sonia Foschiatti\\
	Faculty of Mathematics\\
	University of Vienna\\
	Oskar-Morgenstern-Platz 1\\
	A-1090 Vienna, Austria\\
        sonia.foschiatti@univie.ac.at
        \and
	Axel Kittenberger\\
	Faculty of Mathematics\\
	University of Vienna\\
	Oskar-Morgenstern-Platz 1\\
	A-1090 Vienna, Austria\\
        axel.kittenberger@univie.ac.at
}

\keywords{virtual unwrapping; X‑ray computed tomography; educational laboratory; applied mathematics; image reconstruction; genetic algorithms; splines.}

\begin{document}

\makePtitlepage
\makePtitle

\begin{abstract}
    In this paper, we describe ``Archeolab: Computed Tomography and Virtual Unwrapping of Scrolls'', an interdisciplinary undergraduate course guiding students through the virtual reading of scrolls. By studying the theoretical foundations of X-ray Computed Tomography, genetic algorithms, and conducting hands-on experiments and programming sessions, we aim to encourage students to study applied mathematics and promote active participation in lectures. 
\end{abstract}

\listkeywords

\section{INTRODUCTION}\label{sec:intro}

Applied mathematics is attracting increasing attention in undergraduate and graduate curricula (see \cite{Nis87}). The term 'applied' means that the mathematical problem under study, or, more precisely, the situation and questions that define it, belong to the real world, i.e., anything that does not directly involve the domain of mathematics (\cite{Pol79,BluNis91}). An applied mathematical problem requires a real-life scenario to be simplified and translated into mathematical terms, namely a modeling process. However, it has been pointed out that the modeling practice does not play a significant role in undergraduate education, in part because both students and teachers face difficulties in implementing it (\cite{BluFer09}). Therefore, there is a need for courses that overcome these obstacles by including interdisciplinary content and modeling.

In this article, we describe an application-inspired course called 'ArcheoLab: Computed Tomography and Virtual Unwrapping of Scrolls', which has been proposed for a Master's level mathematics curriculum at the University of Vienna. We will also outline how it can be adapted for a Bachelor's level course. Taking inspiration from the archaeological challenge of unraveling the mysterious content of ancient papyri using modern technology and the recent advances in computer science and computer vision, the course aims to motivate students to engage with the study of applied mathematics. More specifically, the selected example was the Vesuvius Challenge \cite{Wik}, a public competition involving the use of X-ray scanning data and machine learning algorithms aiming at developing an automated software pipeline for reading the Herculaneum papyri, carbonized scrolls that were discovered near Naples in the 18th century, virtually. In \cite{KitMinSch22} and \cite{FosKitSch25}, an educational laboratory was outlined to teach X-ray computed tomography via the archaeological application. This was the starting point for the development of our course, which not only implemented a practical experiment with the students but also integrated the theoretical and programming knowledge to build a complete reconstruction pipeline, from data acquisition to virtual unwrapping. Students engage with theoretical results and modeling decisions while developing independent problem-solving strategies, abstract and critical thinking, and computational skills.

The University of Vienna categorizes its curriculum into courses with non-continuous assessment, such as traditional lectures concluding with a single exam, and courses with continuous assessment. ArcheoLab was delivered in the framework of a Master's seminar, a format with continuous assessment that promotes students' self-study and autonomous research on a given topic. To move beyond the standard read-and-present format, we made the course more interactive by combining instructor-led mini-lectures with a laboratory part. Assessment emphasized both research communication and self-study, as well as the practical combination of theoretical and coding exercises.

To support reproducibility and reuse by other instructors, we provide a detailed timeline of the course, lecture slides, datasets, example codes, laboratory instructions, results and feedback from the students, and suggestions for further improvements. We recommend that the instructors have a solid knowledge of inverse problems and applied mathematics, as well as coding competencies in Python or Matlab. For Bachelor 's-level adaptation, we recommend reducing the theoretical prerequisites. Moreover, in the case that the laboratory experience is difficult to implement, we offer a series of databases that can be used by the instructors together with the full reconstruction pipeline.

To outline the organization, we describe the three main didactical blocks.
\begin{enumerate}
	\item  In the first part, each weekly lecture was divided into two sessions. The first covered theoretical topics related to the ArcheoLab: X-ray computed tomography, inverse problems and ill-posedness, evolutionary algorithms, and splines. The second session was a programming laboratory, where we assigned weekly notebooks, combining theoretical questions with coding exercises. 
	\item The second part was dedicated to handicraft, laboratory visit, and virtual reconstruction using our code. Students were divided into five groups, and each group reconstructed a different scroll. Since working directly with X-ray data was infeasible, we simulated carbonized scrolls by wrapping transparent films printed with text.
	\item In the third part, the students were the main actors. Each group selected a topic from a given list of twelve options and gave a presentation in the classroom connecting the chosen theoretical topic and the experimental results obtained.
\end{enumerate}

From a pedagogical point of view, we adopted a flipped classroom approach, a well-known tool that helps increase student interaction during the lecture, and real-time feedback (see \cite{ChuHewChe17}). Before class, students were assigned videos or parts of book chapters, along with questions, to introduce the fundamentals of a given topic. The content was partially recalled during the lecture, but more time was dedicated to advanced topics and in-class group activities. To further promote active learning, we have introduced multiple-choice quizzes and open questions. On the computational side, students have used programming tools and packages in Python. The students were then assigned weekly programming tasks to consolidate their understanding of the code and concepts, and progressively prepare them for the reconstruction pipeline.

The paper is organized as follows. In Section \ref{sec:class} we detail the course design and journey. Subsection \ref{subsec:part1} presents the theoretical and programming preliminaries, including the mathematical modelling of X-ray Computed Tomography and an introduction to evolutionary algorithms. Subsection \ref{subsec:laboratory} describes the hands-on laboratory and the ArcheoLab pipeline from data acquisition to virtual unwrapping. Section \ref{sec:methods} outlines the instructional methods, technologies, and materials, together with links to online resources. Section \ref{sec:results} summarizes students' presentations and projects and reflects on feedback from both the instructors and students.

\section{STRUCTURE OF THE COURSE}\label{sec:class}
We offered the course "Archeolab: Computed Tomography and Virtual Unwrapping of Scrolls" as a seminar during the summer term 2026 at the University of Vienna. We met once per week throughout the semester for a total of 11 sessions of 90 minutes each. Successful completion of the course was credited with 4 ECTS points, totaling 100 hours of work, including time spent in class. In accordance with the University of Vienna's course format definitions, a seminar employs continuous assessment and requires participants to engage independently with a subject area, demonstrating their acquired knowledge through oral and written contributions such as presentations and seminar papers.

All 14 enrolled students were master's students in Mathematics, taking the course as an elective. Although it was offered as a master 's-level seminar due to administrative constraints, its content is equally suitable for advanced undergraduate students. The prerequisites comprise calculus, linear algebra, functional analysis, and foundational programming experience in Python, which suggests the course would also serve well as a third-year elective within the Bachelor's program in Mathematics at the University of Vienna.

\begin{figure}
	\centering
	\includegraphics[width=0.5\linewidth]{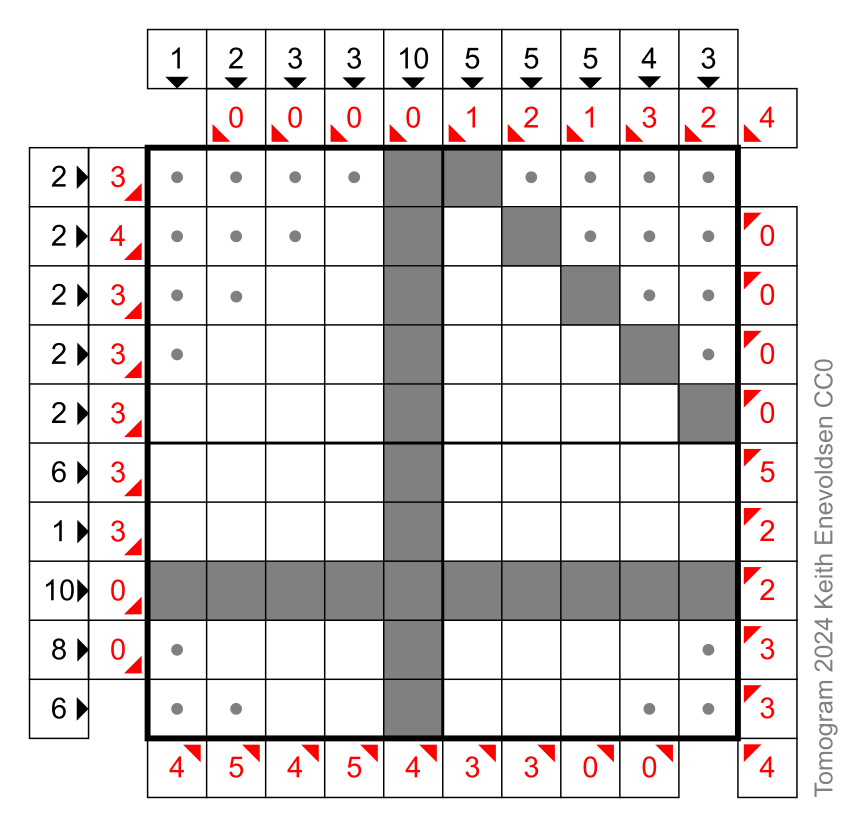}
	\caption{Partially solved tomogram, courtesy of \cite{wlonk_tomograms}.}
	\label{fig:tomogram}
\end{figure}

We opened the course with an introduction to the Vesuvius Challenge and the problem of virtually reading carbonized scrolls from CT data, aiming to immerse students in an active, real-world research problem with open challenges yet to be resolved. To immediately engage the students in the underlying concept of an inverse problem while establishing a collaborative environment, we transitioned from this introduction into a group icebreaker activity. Students were divided into arbitrary groups and tasked with solving "tomogram puzzles" \cite{wlonk_tomograms}. These puzzles act as a simplified, discrete 2D analog of computed tomography, where the goal is to reconstruct a grid of binary values based only on their row, column, and diagonal sums, as illustrated in \autoref{fig:tomogram}. This provided an accessible, gamified entry point into the mathematics of projection and reconstruction, while allowing the students to get to know each other. 

We then presented a concrete, accessible experiment introduced in \cite{FosKitSch25} to replicate the core challenges of the Vesuvius Challenge using simple resources. In that work, a transparent scroll is illuminated with visible light from a projector and recorded with a camera. It was demonstrated that tools from computed tomography, such as filtered backprojection, yield a reasonable reconstruction of the scroll, even when scattering effects of light are not accounted for.

\begin{figure}[htbp]
	\centering
	\includegraphics[width=0.9\textwidth]{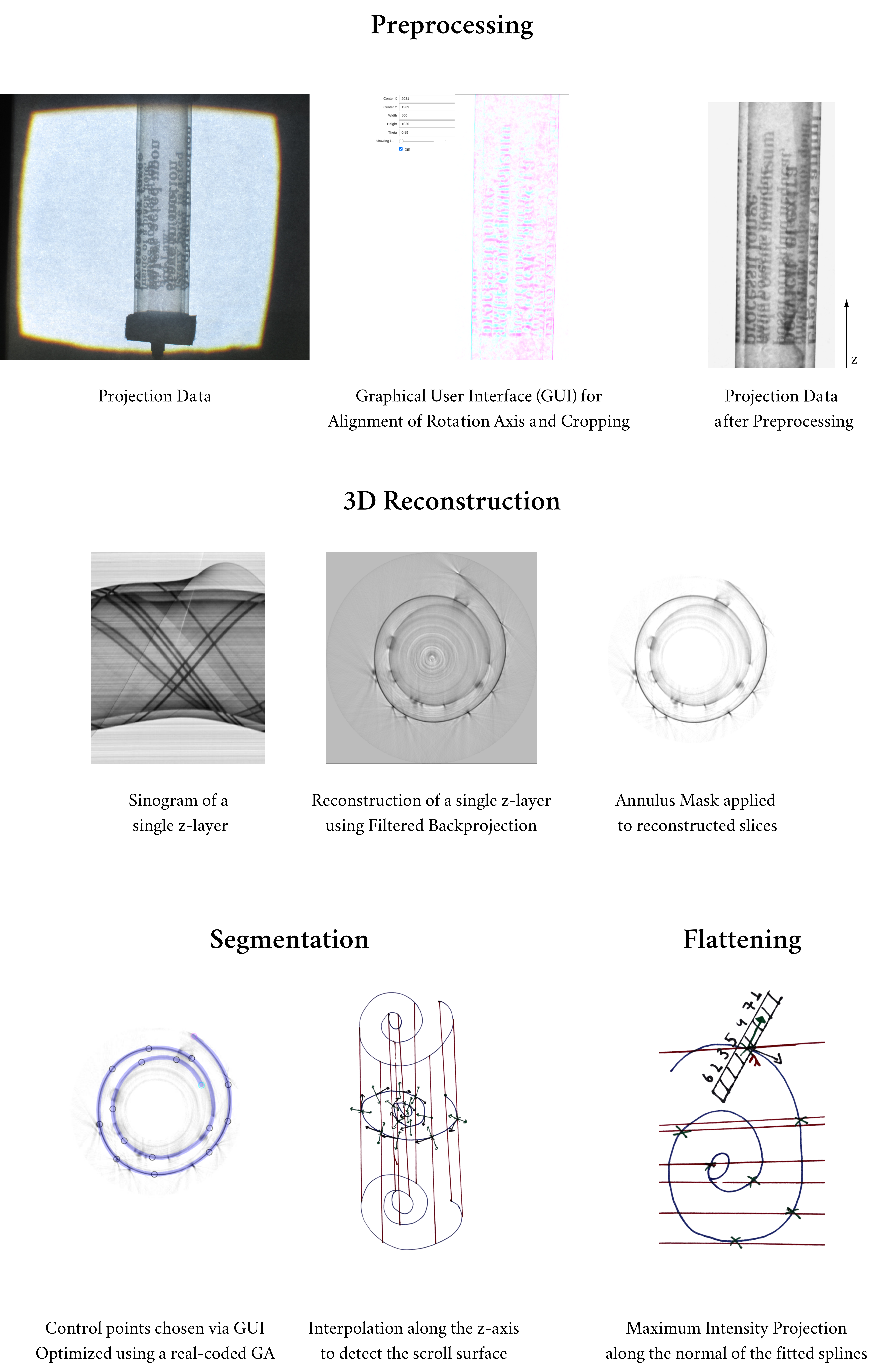}
	\caption{The standard ArcheoLab virtual unwrapping pipeline. The process begins with preprocessing the CT projection images, followed by 3D reconstruction using filtered backprojection. The scroll layers are then segmented and flattened. Pictures of the sinogram, as well as the hand-drawn sketches, are courtesy of our students Aron Riemenschneider, Adrian Killian Raum, Odysseas Konstantatos.}
	\label{fig:pipeline}
\end{figure}

Building on this motivation and the intuition gained from the tomogram activity, we then 
provided an overview of the full reconstruction pipeline for our experiment, which we adapted from \cite{FosKitSch25}. The goal of 
the ArcheoLab pipeline is to virtually read the text inscribed on a physically inaccessible 
scroll. Starting from a set of projection images of the rolled scroll taken from different 
angles, the pipeline proceeds in four stages, summarized in \autoref{tab:reconstruction_steps} and illustrated in \autoref{fig:pipeline}.

\begin{table}[ht]
	\centering
	\begin{tabular}{p{0.25\textwidth}p{0.35\textwidth}p{0.30\textwidth}}
		\toprule
		\textbf{Pipeline stage} & \textbf{Description} & \textbf{Mathematical topics} \\
		\midrule
		Preprocessing & The acquired projection images are loaded and aligned in preparation for 3D reconstruction. & Noise modeling, filter design \\
		\midrule
		3D reconstruction & Filtered backprojection is applied slice by slice to reconstruct the 3D volume of the scroll. & Radon transform, filtered backprojection \\
		\midrule
		Segmentation & For two or more layers, the shape of the scroll in a single layer is interpolated with a B-spline; a real-coded genetic algorithm is then used to optimize this interpolation. & Splines, evolutionary algorithms \\
		\midrule
		Flattening & The scroll surface is detected via interpolation between segmented slices and flattened using a maximum intensity projection. & Splines, maximum intensity projection \\
		\bottomrule
	\end{tabular}
	\caption{Overview of ArcheoLab pipeline stages, their descriptions, and corresponding mathematical topics.}
	\label{tab:reconstruction_steps}
\end{table}

The course was divided into three parts: a theory and programming block covering CT reconstruction and evolutionary algorithms, a hands-on scroll laboratory, and student presentations, each described in the subsections below. In \autoref{tab:course_overview}, we present an overview of the timeline and the used resources. All course materials, including lecture notes, Python notebooks, and problem sets, as well as instructor resources, are freely available on the course homepage, see \autoref{section:material}.

\begin{table}[ht]
	\centering
	\begin{tabular}{p{1cm}p{6cm}p{6cm}}
		\toprule
		\textbf{Session} & \textbf{In-class activities} & \textbf{At-home activities} \\
		\midrule
		1 & Introduction to Vesuvius Challenge, organizational information, tomogram group task & Problem Set 1 \\
		\midrule
		2 & Theory: Lecture 1, Programming: Python Tutorial 1 & Problem Set 2 \\
		\midrule
		3 & Theory: Lecture 2, Programming: Python Tutorial 2 & Problem Set 3 \\
		\midrule
		4 & Theory: Lecture 3, Programming: Python Tutorial 3 & Problem Set 4 \\
		\midrule
		5 & Theory: Lecture 4, Programming: Python Tutorial 4 & Problem Set 5 \\
		\midrule
		6 & Lab visit, Programming: Python Tutorial 5, full ArcheoLab pipeline & Individual work on recorded data and presentations \\
		\midrule
		7 & Q\&A / support with presentations, List of theoretical topics & Group work on recorded data and presentations \\
		\midrule
		8 & Q\&A / support with presentations & Group work on recorded data and presentations \\
		\midrule
		9 & Presentations: Groups A and B & Group work on recorded data and presentations \\
		\midrule
		10 & Presentations: Groups C, D, and E & Feedback form \\
		\midrule
		11 & Wrap-up \& discussion of feedback & --- \\
		\bottomrule
	\end{tabular}
	\caption{Overview of in-class and at-home activities across the 11 course sessions.}
	\label{tab:course_overview}
\end{table}

\subsection{Part 1: Theoretical and programming preliminaries}\label{subsec:part1}
Building on the introductory overview, the first major block of the course focused on the foundational mathematics and computational tools required for the ArcheoLab pipeline. Rather than a traditional lecture format, we employed a flipped classroom approach where students engaged with preparatory materials before class. This allowed in-class sessions to focus on interactive programming tutorials and collaborative problem solving. The goal was to enable the students to later adjust the pipeline in an informed way. 

\subsubsection{The mathematical modeling of X-ray Computed Tomography}
\textit{Material: Sessions 2--3; Lectures 1--2; Python Tutorials 1--2; Problem Sets 1--3.}

This part covered the theory of X-ray CT and its formulation as an inverse problem. Students explored the Radon transform, the Fourier Slice Theorem, Filtered Back-Projection (FBP), and the ill-posedness of the reconstruction problem. Through hands-on Python tutorials, students implemented these concepts to observe the visual effects of sparse angle sampling, limited angle tomography, and Additive White Gaussian Noise. Practical coding exercises emphasized the necessity of regularization. Students systematically compared different filters, and we guided them through extending 2D FBP algorithms to 3D volume reconstructions.

\subsubsection{Evolutionary Algorithms}
\textit{Material: Sessions 4--5; Lectures 3--4; Python Tutorials 3--5; Problem Sets 3--5.}

The second major block of the course introduced evolutionary algorithms and B-splines to optimize the segmentation of the scroll layers. Students progressed from implementing a Simple Genetic Algorithm from scratch to utilizing the DEAP framework \cite{ForRaiGarParGag12} for real-coded algorithms. During this process, they explored various selection strategies, such as Roulette Wheel and Tournament Selection. The theoretical discussions included an elaborate examination of different types and choices of evolutionary algorithms, alongside Holland’s Schema Theorem and the drawbacks of binary encoding. Finally, the module concluded with the introduction of B-splines. Students learned how to mathematically represent and trace the scroll layers identified by the genetic algorithm to compute the normal vectors required for the virtual unwrapping process.

\subsection{Part 2: Hands-on scroll laboratory and virtual reconstruction} \label{subsec:laboratory}
\textit{Material: Full Archeolab Pipeline.}

Building on the theoretical foundations, the second part of the course focused on the practical application of the material learned. The core experimental design utilized in this laboratory, including the low-cost experimental setup involving transparent films printed with text and rolled into a scroll-like object as well as the Python code, was originally developed and detailed in \cite{FosKitSch25}. Providing these materials enabled students to experience the complete virtual unwrapping workflow in a controlled and accessible environment. Students were organized into small groups, each working with their own unique scroll data. The primary objective was not merely to produce a readable flattened image, but rather to analyze the interactions between different pipeline stages and to identify potential sources of error and areas for improvement.

In preparation for this phase, students formed groups of two to three people during Session 4 and selected the texts for their respective scrolls. We printed these texts onto transparent film and provided the necessary materials for students to physically construct their scrolls. To optimize classroom time and accommodate logistical constraints, the instructors pre-recorded the projection data while the students concluded their theoretical coursework. This decision allowed the course to remain strictly focused on the mathematical and computational aspects of the reconstruction. However, the curriculum could be adapted to include the data acquisition process as an active student task.

We initiated the practical component in Session 6 with a visit to the tomography lab, where the data acquisition was performed. Students had the opportunity to examine the experimental setup, take photographs, and ask questions. It is worth noting that the experiment does not require a specialized, traditional laboratory. It relies on accessible equipment and can be easily conducted in a standard office space \cite{FosKitSch25}.

After a comprehensive overview of the entire pipeline, each group was assigned the projection data for a scroll constructed by another group. This cross-assignment introduced an element of authentic discovery. Because the students were unaware of the hidden text, they were challenged to genuinely decipher it. Equipped with the complete dataset and a functional algorithmic framework, the students dedicated the following two sessions, as well as independent study time, to their reconstructions, with instructors available in class to provide targeted guidance and troubleshoot technical issues.

\subsection{Part 3: Students' presentations}\label{subsec:presentations}
\textit{Material: List of theoretical topics. }

We concluded the course with final presentations over two weeks, in which each group showcased both their practical reconstruction results and their independent research on a theoretical topic. 

For the practical component, students presented the deciphered text from their assigned scroll and discussed the specific challenges they encountered, explaining their algorithmic parameter choices and presenting any extra preprocessing or filtering steps they added to the pipeline to enhance the final results.

For the theoretical component, each group was required to choose a topic of choice or select a topic from a curated list of eleven options provided by the instructors. This list was designed to extend the foundational material covered in the lectures and included a mix of mathematical and applied subjects. By combining their experimental findings with these advanced theoretical concepts, the students were able to contextualize their practical work within the broader scope of applied mathematics.

A key pedagogical goal was to foster a supportive and safe environment during these final sessions. We emphasized that encountering errors is an inherent part of the research process, and we actively encouraged students to share approaches that did not yield the expected results. In applied mathematics, particularly when handling real-world data, every failed reconstruction or visible artifact presents a valuable learning opportunity.

\section{IMPLEMENTATION: METHODS AND RESOURCES}\label{sec:methods}

In this section we give additional details about the instructional methods and tools used during the course, complementing the course structure described in Section \ref{sec:class}.

\subsection{Instructional methods}

From a didactical perspective, we adopted a flipped-classroom approach. From the beginning, students engaged with short preparatory materials related to the theoretical topics (videos of length 3-10 minutes and brief readings) and were asked to solve a few targeted exercises before class. 
Class time was then dedicated to small-group problem solving (teams of 2-4 students), class discussions, and short quizzes.
Quizzes were provided using different platforms, like Slido and Kahoot. During the first lectures, quizzes addressed different questions related to why the students decided to attend this course and what do they expect from it, and preparatory questions about the theoretical content and the programming content. After the completion of the theoretical part, the instructors proposed a final quiz helping them assessing the overall understanding of the theory and the main concepts of the course.

To simulate authentic research and avoid a purely evaluative classroom, we established a supportive environment that normalized error as part of the learning and research process. In grading the assignments, we focused on the completion of the tasks and we emphasized a constructive and process-oriented feedback, giving suggestions for improvement in the resolution of the task or appreciating an original solution. During the reconstruction phase and the oral presentations we actively encouraged students in solving their tasks, underline the inherent difficulty of the overall process. In the final lecture, we asked students to share approaches or parameter choices that did or did not yield the expected results. This is inline with applied mathematics practice: handling noisy real-world data, failed reconstructions, algorithm crashes, or visible artifacts allow to get theoretical and practical insights, becoming a valuable learning resource.

To evaluate the oral presentations of the students, the instructors used an evaluation matrix based on the British Council sample oral presentation marking criteria, scored on a 1-5 scale. The seven evaluation criteria were: delivery (clarity, audibility, pacing and engagement with the audience); knowledge and content (evaluation of confidence in the material, evidence of personal research, ability to explain key ideas); structure (logical flow, clarity of arguments, coverage of the relevant points); analysis and evaluation (quality of analysis, interpretation of the material, acknowledgment of limitations); use of visual aids (clarity and appropriateness of figures and animations); response to questions (accuracy and appropriateness when addressing peers' and instructors' questions); and teamwork (coherence of the group presentations and integration of individual contributions). To promote peer learning and feedback during the oral presentations, students received a feedback matrix that guided constructive, peer-to-peer comments.

\subsection{Technologies and materials}

Course materials were provided via Moodle to share slides, code templates, datasets, and references; to collect homework, group choice, argument choice and feedback submissions; and to keep an active communication with the students. 

The programming language adopted for the ArcheoLab was Python beacuse it is open-source, easy to learn, and supported by the libraries needed for the course's theoretical and applied components. Morevoer, an extensive online documentation supports beginners and enables interested students to explore advanced libraries. We remark here that programming materials can be easily adapted to MATLAB. However, this programming language requires licensed access to relevant toolboxes and attention to differences in defaults values and numerical behavior, making Python a more democratic choice.

Notebooks were presented through the University of Vienna's JupyterHub, giving students browser-based access to a preconfigured environment without local installation. This allowed instructors and students to work in real time on the platform. During class, students worked through short, scaffolded tasks; more demanding programming exercises were assigned for homework. Across the semester, we used a total of five lecture notebooks, one additional homework notebook, and the ArcheoLab pipeline notebook. For the lecture notebooks we relied on scikit-image for tomographic operations and image processing (in particular the functions \texttt{radon} and \texttt{iradon} for the Radon transform, the Filtered Back-Projection), DEAP for genetic algorithms, and standard scientific libraries like NumPy, SciPy and Matplotlib. The Archeolab notebook required a more extensive set of dependencies; complete environmental details, and package versions are provided on the course webpage.

\subsection{Useful links to material and resources}\label{section:material}

\begin{itemize}
	\item Instructor resources (organizational slides, theoretical notes, Jupyter Notebooks, and bibliography): \url{https://csc1.gitlab.io/archeolab/archeolab-course/}.
	\item Laboratory setup description: \url{https://csc1.gitlab.io/archeolab/lab/}.
	\item Datasets (Latin and Greek scrolls): \url{https://gitlab.com/csc1/archeolab/-/tree/master/scrolls/Recordings?ref_type=heads}
	\item Latin scroll tutorial: \url{https://csc1.gitlab.io/archeolab/latin/unwrapping/}.
\end{itemize}

\section{RESULTS AND FEEDBACK}\label{sec:results}

\subsection{Reconstruction Results}
During the final project phase, all student groups successfully ran the ArcheoLab pipeline to virtually unwrap their assigned scrolls. All groups used JupyterLab with our custom kernel for their reconstructions, with one group additionally working with a local installation. 

One important lesson for us as instructors, as well as for the students, was that the fine, stationary texture of the paper screen led to severe ring artifacts in the reconstruction. Scaling the projection data by a factor of 0.5 averaged adjacent pixels together, acting as a low-pass filter that significantly improved the final reconstructions.
Beyond this shared preprocessing step, students optimized various parameters to suit their datasets. As shown in \autoref{fig:student-reconstructions} (B) and (D), two groups reduced the number of generations in the genetic algorithm significantly from 100 to 20. Another group investigated the effects of enhancing the recorded images by applying Gaussian smoothing and sharpening prior to the 3D reconstruction (\autoref{fig:student-reconstructions} (C)).

One group performed custom destriping on their sinograms based on \cite{VoAtwDra18}, combined with additional two-dimensional smoothing to further improve the aforementioned ring artifacts (\autoref{fig:student-reconstructions} (A)). This was also the only team to substitute the default Ramp filter with a Shepp-Logan filter, and to increase interpolation accuracy by segmenting three separate layers instead of the standard two.
Finally, another group focused heavily on post-processing. Utilizing the knowledge gained from their theoretical topic on image segmentation, they applied a variety of algorithms to the unwrapped scroll, including Canny Edge detection, Otsu's method, and Region Growing, to isolate the ink and maximize text legibility (\autoref{fig:student-reconstructions} (E) and (F)).

\begin{figure}[htbp]
	\centering
	\includegraphics[width=\textwidth]{reconstructions.pdf}
	\caption{Summary of the virtual unwrapping results achieved by the student groups.}
	\label{fig:student-reconstructions}
\end{figure}

\subsection{Instructor and Student Feedback}

\subsubsection{Origins and designs}

This course was novel for both instructors and posed several challenges in its design, development, and assessment. The idea of offering an applied mathematical course originated with the first author after being invited to collaborate with the second author on a student laboratory component for another seminar. With the interest and support of two faculty members, the authors were allowed to develop and teach the course independently. Planning took several months and focused on content selection, course format, and appropriate didactical methods. From the beginning, the instructors agreed to divide the course into four components: a theoretical part, a programming part, a hands-on laboratory, and a final assessment. Due to time and contractual constraints, the final assessment consisted of a presentation integrating both theoretical and experimental aspects. 

\subsubsection{Instructor reflections}
Overall, the instructors' impressions were very positive. Students engaged actively in class discussions, completed assignments during the semester on time, and generally appreciated the opportunity to work in groups. The balance between theory and practice appeared to motivate participation and support learning. 

During the semester, the instructors faced several challenges. In terms of content and pacing, they needed to calibrate the workload across the theoretical and programming parts and ensure that students with varied backgrounds could follow the lectures and study the material independently. Since the instructors did not know students' prior curricula, they avoided going too deeply into some theoretical topics and provided more scaffolding in the programming part. Another challenge was posed by the choice of bibliography: while the instructors had strong expertise for the first theoretical topic aligned with their research group's interests, identifying high-quality resources for the genetic algorithm part required additional effort. 

During the course, we had to cover a large number of advanced topics, which occasionally felt like too many given the limited amount of time. One concrete option to reduce the overall workload is to simplify the segmentation phase. We realized when working with the student datasets that running the genetic algorithm to optimize the segmentation is not strictly necessary. The reconstruction also worked effectively if the manual control points were chosen with care. This approach would allow instructors to set aside more time to discuss splines in greater depth.

For other instructors planning a similar course, we recommend careful advance planning of prerequisites and resources. The course benefits from a background in inverse problems, a solid confidence in programming, and early planning of the experimental activities.

\subsubsection{Student feedback}

Fourteen students enrolled and twelve successfully completed the course. Eleven completed the standard university evaluation (response rate 78.6\%) and twelve completed an additional survey elaborated by the instructors consisting of 15 questions (6 multiple choice and 9 open questions). The optional university course evaluation, which includes a program-level section and a course-quality section, reported that 72.3 \% of students rated the course "very good", and 27.7\% rated it "good". Regarding the workload/ ECTS, 81.8 \% found the allocation adequate, while the remainder judged it slightly low for the required work. Concerning teaching quality, most students agreed that the instructors explained difficult concepts clearly, fostered discussion and engagement, structured the content well, and maintained a positive classroom atmosphere. 

The instructor survey complemented the standard university one by probing motivation for course selection, alignment with students' interests, perceived confidence in what was learned and experiences with group work. 

75\% of the students took the course in the second semester of the Master's program, while 25\% took it in the fourth semester. The main reason for taking the course was for personal interest (58.3\%) and elective fitting with the curriculum (41.7\%). The majority of the students choose the course in connection to personal interests or future goals (10/12).

Students express appreciation for the hands-on activities (uncommon in other courses), the collaborative and interactive nature of the course, the variety of tasks in the assignments, the coding component, group activities, the flipped classroom approach, the independent project, and the friendly class atmosphere. 

Students generally valued group work for peer learning and problem-solving. The majority  perceived that the work distribution was rather fair (58.3\%), several students percieved it fair (25.0\%) and a minority neither fair nor unfair (8.3\%) and rather unfair (8.3\%). 58.3\% of the students estimated the share of work done personally around 40-60\% of total group work, whereas 41.7\% about 25-40\% . Several students enjoyed the final group work, because it allowed them to get deeper into the understanding of a topic, and leading to a successful project.

Concerning the relevance of the course to personal interests or future goals, students shared different opinions. For some, the course served as an introduction to tomography and inverse problems and sparked new interest in the subject. Other valued the programming and project component positively, as well as the demonstration of how mathematics applies across domains, which, in some cases, inspired plans for a Master's thesis in tomography or medical imaging. A few students noted that, although the topic was not directly related to their future goal, they still found the course tasks motivating and valuable.

Overall, about half of the students reported that they can follow and explain the main methods and results covered in the course, several indicated that they can discuss the methods in detail, and a smaller group felt able to critically evaluate and extend the methods on its own.

Finally, the students proposed to keep the flipped format of the lectures, the small-group, collaborative environment and group projects, the use of JupyterHub and Jupyter notebooks, and the hands-on assignments. What students appreciated most of this course was the applied focus of the course, the end-to-end, real-data pipeline experience and the support and enthusiasm of the instructors.   

\section{DISCLOSURE STATEMENT}

During the preparation of this article and the course materials, the authors used assistive and generative artificial intelligence tools (Grammarly, Perplexity (powered by Gemini 3.1 Pro and Claude Sonnet 4.6 Thinking), AI of the University of Vienna).

AI tools were strictly utilized for editorial, technical, and administrative assistance. No AI tools were used for idea generation, conceptualization, or conducting the underlying research. For manuscript preparation, the aforementioned tools were used to assist with language editing, improving readability, structuring paragraphs, LaTeX formatting, and summarizing the course structure. For the development of the pedagogical materials discussed in the paper, AI was used to generate reading assignments, debug code, and construct coding exercises.

After using these tools, the authors carefully reviewed, edited, and verified all generated text and code. The authors take full responsibility for the originality, validity, accuracy, and integrity of the content presented in this article.

No potential conflict of interest was reported by the authors.

\section*{ACKNOWLEDGMENTS}

This research was funded in whole or in part by the Austrian Science Fund (FWF) SFB 10.55776/F68: "Tomography Across the Scales", projects F6804-N36 (Quantitative Coupled Physics
Imaging) and F68-01 (Coordination Project). For open access purposes, the author has applied a CC BY public copyright license to any author-accepted manuscript version arising from this submission. The authors would like to thank all the students that actively attended the course for their contributions and valuable insights.

\section*{Biographical Sketches}

\vspace*{.3 true cm}
\noindent \textbf{Denise Schmutz} is a Doctoral Candidate at the Faculty of Mathematics, University of Vienna, specializing in motion reconstruction in tomographic imaging. She is passionate about science outreach, teaching the mathematics of computed tomography at the Vienna Children's University. In 2026, her research was featured in the FÄKT program of the Austrian Academy of Sciences, a science communication initiative aimed at young audiences.

\textbf{Sonia Foschiatti} is currently a Postdoctoral Researcher at the Faculty of Mathematics, University of Vienna, Austria. Her research interests include inverse problems for partial differential equations, and outreach activities for the topics concerning the Special Research Program ``Tomography across the scales'', like the Archeolab.

\vspace*{.3 true cm} 
\noindent \textbf{Axel Kittenberger} is currently an IT specialist at the Faculty of Mathematics, University of Vienna, Austria.

\end{document}